\title{Improved Pebbling Bounds}
\author{\small Melody Chan \\
\small Yale University \\
\small \tt{melody.chan@aya.yale.edu} \and
\small Anant P.~Godbole \\
\small Department of Mathematics \\
\small East Tennessee State University \\
\small \tt{godbolea@etsu.edu}}

\documentclass [12 pt]{article}
\usepackage{amsfonts}
\begin{document}
\def\diam{{\rm diam}}
\def\dist{{\rm dist}}
\def\lr{\left(}
\def\rr{\right)}
\def\lc{\left\{}
\def\rc{\right\}}
\def\qed{\vbox{\hrule\hbox{\vrule\kern3pt\vbox{\kern6pt}\kern3pt\vrule}\hrule}}
\newtheorem{thm}{Theorem}
\newtheorem{result}{Result}
\maketitle
\begin{abstract}
Consider a configuration of pebbles distributed on the vertices of a connected graph of order $n$.  A pebbling step consists of removing two pebbles from a given vertex and placing one pebble on an adjacent vertex.  A distribution of pebbles on a graph is called solvable if it is possible to place a pebble on any given vertex using a sequence of pebbling steps.  The pebbling number of a graph, denoted $f(G)$, is the minimal number of pebbles such that every configuration of $f(G)$ pebbles on $G$ is solvable.  We derive several general upper bounds on the pebbling number, improving previous results.
\end{abstract}

\section{Introduction}
	\subsection{Definitions}

	Given a connected graph $G = (V,E)$, let $D: V \rightarrow \mathbb{N}$ be a distribution of $\sum_v D(v)$ identical pebbles on the vertices of $G$.  A {\em pebbling step} consists of removing two pebbles from a given vertex and placing one of these pebbles on an adjacent vertex (the other pebble is removed from the graph.)  Given a root vertex $v \in V$, we say that $D$ is {\em $v$-solvable} if we can place at least one pebble on $v$ after some number of pebbling steps.  $D$ is called solvable if $D$ is $v$-solvable for all $v \in V$.  Let the {\em size} of a distribution $D$ be $|D| = \sum_v D(v)$, the total number of pebbles on the graph.  Then the pebbling number $f(G)$ is defined to be the smallest integer $N$ such that any distribution of size $N$ is solvable.
Also, we define $f(G,v)$ to be the smallest integer $N$ such that any distribution of size $N$ is $v$-solvable.
	Throughout this paper, $G = (V,E)$ denotes a simple, connected graph, $n = |V|$ is the number of vertices in $G$, and $d = \diam(G)$ is the diameter of $G$.

	\subsection{Known Bounds on the Pebbling Number}

	First we state an elementary bound on $f(G)$, given in \cite{hurlbert}:  Clearly, if $D(v) = 0$ for the root vertex $v$ and $D(w) = 1$ for all other vertices $w$, then $D$ is unsolvable.  Also, given vertices $v_1$ and $v_2$ at distance $d = \diam(G)$ from each other, if $D(v_1) = 2^d - 1$ and $D(w) = 0$ for every other vertex, then vertex $v_2$ cannot be pebbled and $D$ is unsolvable.  These facts were noted by Chung \cite{chung}.  For an upper bound, note that if $|D| = (n-1)(2^d - 1) + 1$, then either each vertex has at least one pebble on it, or, by the pigeonhole principle, there exists a vertex with $2^d$ pebbles on it, and in either case $D$ is solvable.  To summarize,
\begin{equation}
\max\{n,2^d\} \leq f(G) \leq (n-1)(2^d-1) + 1.
\end{equation}
Note that the upper bound in (1) is sharp if $G$ is the complete graph $K_n$, but is way off target if $G=P_n$, the path on $n$ vertices.

\section{New Upper Bounds}

In this section, we prove four new upper bounds on the pebbling number, two of which always improve the bound stated above, and with the rest doing better than (1) in most cases.

\begin{thm}
\begin{displaymath}
f(G) \leq (n-d)(2^d-1)+1.
\end{displaymath}
\end{thm}
\noindent{\em Proof.}  First, observe that for a path $P_n$ on $n$ vertices, $f(P_n) = 2^{n-1}$.  A simple proof of this result can be found in \cite{hurlbert}.  Now, given any root vertex $v \in V$, consider a set $S_v = \{Q_1, \ldots, Q_m$\} of $m$ paths on $G$, with the following two properties:
\begin{itemize}
\item each path has one endpoint at $v$; and
\item each vertex in $G$ is on at least one path.
\end{itemize}
These paths may overlap, although that is not necessarily the case.  Let $q_i$ be the length of the path $Q_i$.  If some path $Q_i$ has $2^{q_i}$ pebbles on it, then $v$ can be reached using that path.  Thus, the pigeonhole principle guarantees that 
\begin{displaymath} \label{sumoverq}
f(G,v) \leq \left(\sum_{i=1}^{m} (2^{q_i} - 1) \right) + 1.
\end{displaymath}
To construct one such path set, let $Q_1$ be the path from $v$ to some vertex $w$ at maximum distance from $v$, that is, $\dist(v,w) = e(v)$ where $e(v) = \max\{\dist(v,w)\}$ is the eccentricity of $v$.  Then there are $n-e(v)-1$ vertices not on that path.  In the worst case, i.e., in the sense of maximizing the number of paths, each of those vertices requires a distinct path of length no longer than $e(v)$ to connect it to $v$.  Thus, we have in total $n-e(v)$ paths of length at most $e(v)$, and 
\begin{eqnarray*}
f(G,v) & \leq & \left(\sum_{i=1}^{n-e(v)} (2^{e(v)} - 1) \right) + 1 \\
& = & (n-e(v))(2^{e(v)}-1)+1.
\end{eqnarray*}
It follows that
\begin{eqnarray*}
f(G) & = & \max \left\{ f(G,v) \right\} \\
& \leq & \max\left\{(n-e(v))(2^{e(v)}-1)+1 \right\} \\
& = & (n-d)(2^d-1)+1,
\end{eqnarray*}
where the last equality above follows from the fact that the function $(n-j)(2^j-1)$ is monotone increasing for {\it integers} $j$ with $1\le j\le n-1$; to see this, we set $\phi(j)=(n-j)(2^j-1)$ and note that $\phi(j+1)\ge\phi(j)$ if and only if $(n-j)2^j\ge2^{j+1}-1$, and thus if $n-j\ge2$, or $j\le n-2$.

  Note that the upper bound of Theorem 1 obviously improves (1), and is sharp for both $G=K_n$ and $G=P_n$.\hfill\qed

\begin{thm}
\begin{displaymath}
f(G) \leq \left(n + \left \lfloor \frac{n-1}{d} \right \rfloor-1 \right)\left(2^{d-1}\right) - n + 2.
\end{displaymath}
\end{thm}

\noindent{\em Proof.} Given a root vertex $v \in V$ and $k\ge1$, let $\{w_1, \ldots, w_k\} \subset V$ be a set of vertices such that $\dist(w_i, v) = e(v)$ and there exists a set $\{p_1, \ldots, p_k\}$ of length-$e(v)$ paths, with path $p_i$ connecting $w_i$ and $v$, such that no two paths share any vertex except $v$.  Such a set of vertices must exist for, in the worst case we may have $k=1$. Then the number of such paths, $k$, must satisfy $k \leq c$, where $c = \left \lfloor \frac{n-1}{e(v)} \right \rfloor$.    

Now, we claim that 
\begin{displaymath}
f(G,v) \leq (n-1)(2^{e(v)-1} - 1) + c(2^{e(v)-1}) + 1.
\end{displaymath}
To see this, observe that if $v$ already has a pebble on it or if any vertex has $2^{e(v)}$ pebbles on it, then we know that $D$ is $v$-solvable.  If neither is the case, then (with $|D|=(n-1)(2^{e(v)-1} - 1) + c(2^{e(v)-1}) + 1$), the pigeonhole principle guarantees that there exists a set of $c+1$ vertices with at least $2^{e(v)-1}$ pebbles each. This is because we get a contradiction if $v$ has no pebbles, there are at most $c$ vertices with at most $2^{e(v)}-1$ pebbles each, and the remaining vertices have at most $2^{e(v)-1}-1$ pebbles each.  

If any of these $c + 1$ vertices is at distance at most $e(v)-1$ from $v$, we are done.  If all $c+1$ vertices are at distance $e(v)$ from $v$, then there must exist vertices $u_1$ and $u_2$, each with at least $2^{e(v)-1}$ pebbles, with length-$e(v)$ paths $p_{u_1}$ and $p_{u_2}$ to $v$ that both pass through some vertex $w$ adjacent to $v$.  Then we can place 2 pebbles on $w$, and hence we can place 1 pebble on $v$.  Thus, $D$ is $v$-solvable and we have
\begin{eqnarray*}
f(G) & = & \max \left\{ f(G,v) \right\} \\
& \leq & \max \left\{ \left( n + \left \lfloor \frac{n-1}{e(v)} \right \rfloor-1 \right)(2^{e(v)-1}) - n + 2 \right\}\\
& = & \left(n + \left \lfloor \frac{n-1}{d} \right \rfloor-1 \right)\left(2^{d-1}\right) - n + 2,
\end{eqnarray*}
as asserted.\hfill\qed

\bigskip

Note that Theorem 2 is sharp if $G=K_n$, but {\it not} if $G=P_n$; in general, it is easy to verify that Theorem 2 performs better than Theorem 1 whenever $d$ is small and $n$ is large enough.  Also, it is easy to check that Theorem 2 always improves the upper bound given by (1):  to see this, we simply rewrite the latter as $(2n-2)2^{d-1}-n+2$ and note that 
$$(2n-2)2^{d-1}-n+2\ge \left(n + \left \lfloor \frac{n-1}{d} \right \rfloor-1 \right)\left(2^{d-1}\right) - n + 2,$$
with equality holding if and only if $d=1$.

Next, we derive (Theorem 3) a pebbling bound for graphs with {\em efficient dominating sets} which is extended in Theorem 4 to a general bound in terms of the domination number.  Recall that a set $S \subset V$ is said to be a {\em perfect dominating set} if each vertex $v \in V$ is a member of $S$ or is adjacent to exactly one member of $S$.  This definition is not entirely suitable for use in Theorem 3, since we could have edges between vertices in $S$, and $S$ as defined could even be all of $V(G)$.  An {\em independent dominating set} \cite{haynes} is a dominating set whose vertices are independent.  But a vertex $v$ may be adjacent to two or more vertices of an independent  dominating set.  We will thus need the set $S$ to be independent as well as perfect; in other words, we must have an efficient dominating set (see \cite{haynes}, pg.~108, for more on perfect, independent, and efficient domination.)  Our original proof of Theorems 3 and 4 made use of the following result from Pachter et al.~\cite{pachter}:

$$
f(G) \leq n+1 \textrm{ if } \diam(G) = 2,
$$
but, thanks to the suggestions of {\it two} referees, we instead use a fact, proved in \cite{moews}, to prove stronger forms of  Theorems 3 and 4 than in the original version of this paper. Recall that the $k$-pebbling number $f_k(G)$ of a graph $G$ is the minimum number of pebbles that must be placed on $G$ so that any vertex can have $k$ pebbles placed on it in a series of pebbling moves, regardless of the initial configuration of pebbles.  Moews \cite{moews} proved that 
\begin{equation}
f_k(K_{1,{m-1}})=4k+m-3,
\end{equation}
where $K_{1,r}$ denotes the star on $r+1$ vertices.
\begin{thm} Suppose $G$ has an efficient dominating set of size $\gamma$.  Then
\begin{displaymath}
f(G) \leq 2^{d+1}\gamma+n - 4\gamma + 1.
\end{displaymath}
\end{thm}

\noindent{\em Proof}.  Let $S = \{s_1, \ldots, s_\gamma\}$ be an efficient  dominating set.  For each $s_i$, let $A_i = \{s_i\} \cup \{w : w \textrm{ is adjacent to } s_i\}$.  The $A_i$s partition the vertex set $V$, with each $A_i$ having diameter at most 2.  Also, each $A_i$ contains $K_{1,\vert A_i\vert-1}$ as a subgraph.  Thus by (2), if there exists a set  $A_i$ containing at least $2^{d+1}+\vert A_i\vert-3$ pebbles in total, then we may place $2^{d-1}$ pebbles on any vertex in $A_i$.  Furthermore, for $v \notin A_i$, $v$ is at distance at most $d$ from $s_i$, so $v$ is at distance at most $d-1$ from some vertex in $A_i$ -- since $S$ is an efficient dominating set.  Therefore, we may place a pebble on any $v \in V$.  The pigeonhole principle now guarantees that
\begin{eqnarray*}
f(G) & \leq & \left(\sum_{i=1}^\gamma{\left( 2^{d+1}+\vert A_i\vert-4\right)} \right) + 1 \\
& = & 2^{d+1}\gamma+n - 4\gamma + 1,
\end{eqnarray*}
as required.  \hfill\qed

Note that Theorem 3 yields a better bound than that given by (1) if and only if
\begin{equation}
\gamma\le\frac{(2^d-2)n-(2^d-1)}{2^{d+1}-4}=\frac{n}{2}-\frac{2^d-1}{2^{d+1}-4}.
\end{equation}
Now (3) holds unless $d=1$ or $d\ge2$ and $\gamma\sim n/2$. In the former case, (1) {\it does} do  better than Theorem 3 -- but only because the proof of Theorem 3 only used the star structure of the sets $A_i$, and did not exploit the fact that there might be edges between the vertices $w$ adjacent to $s_i$, thus causing the diameter of $A_i$ to conceivably equal one.  Consider the case where $\gamma\sim n/2$.  Since the only connected graphs with domination number $n/2$ are those whose components are the cycle $C_4$ or the corona $G\circ K_1$ (see pp.~41--42 of \cite{haynes} for a further discussion), we see that Theorem 3 does better than (1) ``almost all" the time. 

Our final result shows how Theorem 3 may be generalized, at little cost, to graphs that do not admit a perfect independent dominating set.
\begin{thm}
Suppose $G$ has a dominating set of size $\gamma$.  Then
\begin{displaymath}
f(G) \leq 2^{d+1}\gamma+n - 3\gamma + 1.
\end{displaymath}
\end{thm}

\noindent{\em Proof}.  Let $S = \{s_1, \ldots, s_\gamma\}$ be a dominating set.  For a fixed vertex $v$, pick paths $P_1,P_2,\ldots,P_\gamma$ (of length at most $d$) of the form $P_i=va_1a_2\ldots a_{L_i}w_is_i$. (In case $v$ is a neighbour of $s_i$ or is $s_i$ itself, the path reduces to $vs_i$ or $s_i$, and $w_i$ is defined to be $s_i$ in either case.  If $\dist(v,s_i)=2$, there are no $a_j$s in the path.)  Now define
$$A_1=N_{s_1}\cup \{w_1,s_1\};$$
$$A_2=(N_{s_2}\setminus N_{s_1})\cup \{w_2,s_2\};$$
$$A_3=(N_{s_3}\setminus (N_{s_1}\cup N_{s_2}))\cup \{w_3,s_3\};$$
$$\vdots$$
$$A_i=\lr N_{s_i}\setminus \bigcup_{j=1}^{i-1}N_{s_j}\rr\cup \{w_i,s_i\};$$
$$\vdots$$
$$A_\gamma=\lr N_{s_\gamma}\setminus \bigcup_{j=1}^{ \gamma-1}N_{s_j}\rr\cup \{w_\gamma,s_\gamma\},$$
where $N_x$ consists of all neighbors of $x$ that do not belong to the set $S$.  Notice that the $A_i$s are again sets of diameter at most 2.  The rest of the argument follows the proof of Theorem 3 very closely:  if there exists a set of vertices $A_i$ containing at least $2^{d+1}+\vert A_i\vert-3$ pebbles in total, then we may place $2^{d-1}$ pebbles on any vertex in $A_i$.  Now $v$ is at distance at most $d-1$ from at least one vertex in $A_i$, namely $w_i$.  Therefore, we may place a pebble on $v$.  The pigeonhole principle now guarantees that
\nolinebreak
\begin{eqnarray*}
f(G,v) & \leq & \left(\sum_{i=1}^\gamma{\left( 2^{d+1}+\vert A_i\vert-4\right)} \right) + 1 \\
&=&\gamma2^{d+1}-4\gamma+1+\sum_{i=1}^\gamma \vert A_i\vert\\
&\le&\gamma2^{d+1}-4\gamma+1+
\sum_{i=1}^\gamma\lr\left\vert\lr N_{s_i}\setminus \bigcup_{j=0}^{ i-1}N_{s_j}\rr\cup \{s_i\}\right\vert + 1\rr  \\
& = & \gamma2^{d+1}-3\gamma+1+n,
\end{eqnarray*}
and thus that $f(G)\le 2^{d+1}\gamma+n - 3\gamma + 1$,
as claimed.\hfill\qed

As with Theorem 3, Theorem 4 might not always yield a bound better than that given by (1).  It is easy to verify, however, that this occurs whenever
\begin{equation}\gamma\le\frac{(2^d-2)n-(2^d-1)}{2^{d+1}-3}.\end{equation}
Even for $d=2$, (4) holds if $\gamma\le(2n-3)/5$, which is encouraging since $\gamma\le2n/5$ for ``most" connected graphs with minimum degree at least 2 (see pp.~41--42 of \cite{haynes}).  As $d$ increases, moreover, we see that (4) holds unless $\gamma$ is close to $n/2$.  Thus, Theorem 4 improves on (1) in many cases.    

In a similar fashion, Theorems 3 and 4 may respectively be checked to outperform Theorem 1 if 
\begin{equation}
\gamma\le\frac{(2^d-2)n-d(2^d-1)}{2^{d+1}-4}
\end{equation}
and
\begin{equation}
\gamma\le\frac{(2^d-2)n-d(2^d-1)}{2^{d+1}-3},
\end{equation}
while Theorems 3 and 4 do better than Theorem 2 if, respectively,
\begin{equation}
\gamma\le\frac{(2^{d-1}-2)n-\lr\lfloor\frac{n-1}{d}\rfloor-1\rr2^{d-1}+1}{2^{d+1}-4}
\end{equation}
and 
\begin{equation}
\gamma\le\frac{(2^{d-1}-2)n-\lr\lfloor\frac{n-1}{d}\rfloor-1\rr2^{d-1}+1}{2^{d+1}-3}.
\end{equation}

Theorems 3 or 4 thus often do better than Theorems 1 and 2, but are they ever tight?  This question was raised by one of the referees.  We provide a partial answer.  For the complete graph, our baseline test case, we have $f(K_n)=n$ but using Theorem 3 with $\gamma=1$ and $d=1$ yields a bound of 
$n+2^{d+1}\gamma-4\gamma+1=n+1$.  The next obvious case to check is $K_{1,n}$, whose pebbling number equals $n+2$.  With $d=2$ and $\gamma=1$, however, we see that the upper bound of Theorem 3 is $n+6$, so we only have asymptotic tightness.  The same asymptotic tightness holds if we consider the graph $G$ consisting, for even $n$, of two $K_{1,{n/2}}$'s connected at their roots  by an edge $v_1-v_2$; it is easy to verify that $f(G)= n+6$, with the worst case configuration being no pebbles on $v_1$ or $v_2$, one pebble at each of $(n/2)-1$ vertices on the $K_{1,n/2}$'s and 8 pebbles on any other vertex.  However Theorem 4 applied with $\gamma=2$ yields a bound of $n+29$ for this diameter 3 graph on $n+2$ vertices.  In general we believe that the nature of the proof of either theorem, which uses the pigeonhole principle in a worst case scenario fashion, is unlikely to result in a tight result.  Possible directions for improvement are suggested below.
  
\section{Open Problems}  It would be interesting to develop general bounds on the pebbling number that are not in terms of the diameter of the graph.  We have made some progress in this matter by proving a hybrid bound that depends also on the domination number, but much more needs to be done in this regard.  Of particular interest would be bounds that depend on the girth of the graph, or are expressed in terms of more robust graph invariants such as the tree width (see Robin Thomas' NSF-CBMS lecture notes at {\tt http://www.math.gatech.edu/$\sim$thomas/SLIDE/CBMS/} (a book \cite{thomas} is forthcoming) for an exposition of tree decompositions and tree width).  

Also, how much of an improvement can be made in Theorems 3 and 4 by considering decompositions into sets of (even) diameter four or higher, rather than into the diameter two sets $A_i$ considered in the proofs of these theorems?  This would be necessary if one considers distance $k$-domination, $k\ge2$.  Progress along these lines might be contingent on obtaining tight upper bounds, analogous to those obtained in \cite{pachter} for $d=2$, on the pebbling numbers of graphs with diameter four or higher (the diameter three case has been solved recently by Bukh \cite{bukh}).  Or perhaps we might be able to use the embedded tree structure of the $k$-domination graph of a vertex $v$, together with results in \cite{moews} along the lines of (2), to make the needed improvements.

\section{Acknowledgements}
This work was done under the supervision of Anant Godbole at the East Tennessee State University Research Experience for Undergraduates (REU) site during the summer of 2002, while Melody Chan was a rising sophomore at Yale University.  Support from NSF Grant DMS-0139291 is gratefully acknowledged by both authors.  We would like to thank the three referees whose suggestions for improvement, some of which have been noted at various places in the manuscript, have substantially improved the paper.


\begin{thebibliography}{99}
\bibitem{bukh} B.~Bukh (2005).  ``Maximum pebbling number of graphs of diameter three," to appear in {\it J.~Graph Theory}.
\bibitem{chung} F.~Chung (1989).  ``Pebbling in hypercubes," {\it SIAM J.~Discrete Math.} {\bf 2}, 467--472.
\bibitem{haynes} T.~Haynes, S.~Hedetniemi, and P.~Slater (1998). {\it Fundamentals of Domination in Graphs}, Marcel Dekker, New York.
\bibitem{hurlbert} G.~Hurlbert (1999). ``A survey of graph pebbling," {\it Congress. Numer.} {\bf 139}, 41-64.
\bibitem{moews} D.~Moews (1992).  ``Pebbling graphs," {\it J.~Combin. Theory (B)} {\bf 55}, 244--252.
\bibitem{pachter} L.~Pachter, H.~S.~Snevily, and B.~Voxman (1995). ``On pebbling graphs," {\it Congress. Numer.} {\bf 107}, 65-80.
\bibitem{thomas} R.~Thomas.  {\it Graph Structure and Decomposition.} Forthcoming text.
\end{thebibliography}
\end{document}